\documentclass[fleqn,10pt]{IEEEtran}
\usepackage{multirow}
\usepackage{graphicx}
\usepackage{amssymb}
\usepackage{amsmath}
\usepackage{amsthm}
\usepackage{enumerate}
\usepackage{algorithmic}
\usepackage{algorithm}

\DeclareMathOperator{\sign}{sign}
\DeclareMathOperator{\soft}{soft}

\newcommand{\thh}{\ensuremath{^\text{th}}}
\newcommand{\R}{\ensuremath{\mathbb{R}}}
\newcommand{\xh}{\ensuremath{\hat{x}}}

\theoremstyle{definition}
\newtheorem{defn}{Definition}

\newtheorem{lem}{Lemma}
\newtheorem{prop}{Proposition}

\newenvironment{rcases}
  {\left.\begin{aligned}}
  {\end{aligned}\right\rbrace}

\title{A Penalty Function Promoting Sparsity Within and Across Groups}
\author{\.Ilker Bayram and Sava\c{s}kan Bulek \vspace*{-0.5cm}
\\  \thanks{\.{I}. Bayram is with the Dept. of Electronics and Communications Eng., Istanbul Technical University, Istanbul, Turkey.  E-mail : ibayram@itu.edu.tr. 
S. Bulek is with Qualcomm Atheros, Inc., Auburn Hills, MI, USA. E-mail : sbulek@gmail.com.
}
}
\date{}

\begin{document}

\maketitle

\begin{abstract}
We introduce a new weakly-convex penalty function for signals with a group behavior. The penalty promotes signals with a few number of active groups, where within each group, only a few high magnitude coefficients are active. We derive the threshold function associated with the proposed penalty and study its properties. We discuss how the proposed penalty/threshold function can be useful for signals with isolated non-zeros, such as audio with isolated harmonics along the frequency axis, or reflection functions in exploration seismology where the non-zeros occur on the boundaries of subsoil layers. We demonstrate the use of the proposed penalty/threshold functions in a convex denoising and a non-convex deconvolution formulation. We provide convergent algorithms for both formulations and compare the performance with state-of-the-art methods. 
\end{abstract}

\section{Introduction} \label{sec:intro}

Constraints or prior information derived from sparsity is widely used for regularization in signal processing. Depending on the application domain, the signal of interest may  exhibit additional features than mere sparsity. In this paper, we  consider signals whose coefficients can be clustered in a few  groups where each group itself has few active members. Sparse signals with isolated non-zeros may be considered to fall in that category. We propose a prior function that promotes such signals and demonstrate how to use the function in basic inverse problems of potential interest.

Many natural phenomena can be associated with a sparse underlying process with isolated non-zero components. For instance, the DFT coefficients of a periodic signal are equidistant with respect to the frequency variable. Consequently, quasi-periodic audio signals like speech, music can be represented in the time-frequency domain (via linear transforms \cite{bal13p20}) using components that appear isolated along the frequency axis (i.e., harmonics). 
Another example is related to reflection seismology, where one aims to discover the subsoil layers by sending seismic waves and processing the returning seismic trace \cite{tak12p27}. The seismic trace can be modelled as the convolution of the input seismic wave and the reflection function. The reflection function is a sparse signal with non-zeros occuring due to difference in acoustic impedance between the boundaries of different layers. Since the layers are expected to have some non-zero thickness, the non-zeros, which occur at the boundaries, are isolated. Other than these natural signals, isolated sparsity is also relevant for designed systems. 
For instance, in frequency hopping systems, the parameters of the signal components are constant in between the hopping instances, during which transmission occurs \cite{ang13p64}. Since transmission has to last for some finite amount of time, the hopping instances may be regarded as isolated non-zeros of a sparse signal.

In order to isolate the non-zeros, we work with non-overlapping groups of variables and process the groups independently. We propose a penalty whose threshold function (to be specified below) has the following properties :
\begin{itemize}
\item  If the magnitudes of all variables in a group fall below a threshold, the whole group is set to zero.  
\item Otherwise, a group-dependent threshold is applied so as to eliminate the relatively insignificant coefficients in the group. 
\end{itemize}
The group-dependent threshold serves to separate the large magnitude coefficients from the rest. Specifically, if there are $k$ large-magnitude coefficients in the group, they are kept with little modification, while the rest are set to zero. For $k=1$, this isolates the non-zeros within the group. We remark that this behavior is achieved with a non-adaptive penalty function and without reweighting.

To be more precise, assuming that the size of the groups is $n$, and that $x^{(i)} \in \mathbb{C}^n$ denotes the coefficients belonging to the $i\thh$ group of $x$, we define a penalty function for $\gamma \geq 0$\, as,
\begin{equation}\label{eqn:Pmain}
\mathbf{P}_{\gamma}(x) = \sum_i P_{\gamma}\bigl(x^{(i)} \bigr),
\end{equation}
where for $u \in \mathbb{C}^n$,  $P_{\gamma}$ is defined as,
\begin{equation}\label{eqn:extendN}
P_{\gamma}(u) = \gamma \, \left( \sum_{i=1}^{n-1}\,\sum_{m=i+1}^n |u_i\,u_m| \right) + \|u\|_1.
\end{equation}
The term enclosed in parentheses in \eqref{eqn:extendN} grows rapidly as the number of  large magnitude coefficients in the group increases. Therefore $P_{\gamma}$ strongly penalizes groups containing many large coefficients. 
Given this penalty, we describe how to realize the associated threshold function (or the proximity operator \cite{combettes_chp}) defined for $\lambda \geq 0$ as
\begin{equation}\label{eqn:Tbold}
\mathbf{T}_{\lambda,\gamma}(z) = \arg \min_x \,  \frac{1}{2}\|z -  x\|_2^2 + \lambda\,\mathbf{P}_{\gamma}(x).
\end{equation}
We show that $\mathbf{T}_{\lambda,\gamma}$ is well-defined when ${\lambda\,\gamma < 1}$ and study its behavior. We also show that, as $\gamma \to 1/\lambda$, the threshold function suppresses all but the largest coefficient in each group, provided the magnitude of the largest coefficient exceeds the threshold $\lambda$. We demonstrate the use of the proposed penalty and the threshold function in a  convex formulation for audio denoising and a non-convex formulation for non-blind deconvolution. We provide convergent algorithms for both formulations and demonstrate that  the reconstructions perform favorably compared to those obtained using other penalties/threshold functions.

\subsection*{Related Work}

The proposed penalty function may be regarded as a member of the family of group-based penalty functions (see e.g. \cite{yua06p49,kow09p303,kow09p251,jac09ICML,bay11ICASSP,sie11DAFX,che14p476,rao16p448,sim13p231} for a sample of the literature).
In contrast to our interest, many of these works seek to set whole groups of coefficients to zero, thus achieving sparsity across groups, and do not enforce sparsity within groups. For instance, the $\ell_{2,1}$ norm \cite{kow09p303} is obtained by replacing $P_{\gamma}(x^{(i)})$ with $\|x^{(i)}\|_2$ in \eqref{eqn:Pmain}. The proximity operator associated with the $\ell_{2,1}$ norm sets a whole group to zero if the energy of the group is below a threshold but keeps the group with little modification otherwise. On the other hand, in the Elitist-Lasso (E-Lasso) formulation \cite{kow09p303,kow09p251} (see also \cite{zho10AISTATS} where the method is referred to as Exclusive-Lasso), the target signal contains few non-zeros within each group and sparsity is not enforced across groups. 
Sparsity within groups is also addressed by the sparse-group lasso (SGL) proposed in \cite{sim13p231}. SGL uses a convex combination of an $\ell_1$ norm and an $\ell_{2,1}$ norm as the penalty -- it may also be interpreted as a sum of elastic-net-like penalties \cite{zou05p301}  applied to each group. Therefore SGL uses a convex penalty function. SGL was extended to non-overlapping groups and its performance is thoroughly analyzed in \cite{rao16p448}. 

Another class of related penalties are based on correlations extracted from the observation matrix \cite{tut09p239,anb14p82}.  Given an observation model of the form $y \approx H\,x$, the idea is to derive a positive semi-definite weighting matrix $W$ from the correlations between the columns of $H$\, and use it to define a penalty of the form $x^T\,W\,x$. Since $W$\, does not depend on $x$, the penalty in  \cite{tut09p239,anb14p82} is convex. The targeted effect is uniform treatment of the components of $x$ that produce similar responses. This contrasts with the proposed penalty because if two components of $x$ have similar responses, the proposed penalty $P_{\gamma}$ would prefer to single out one of the components and suppress the other. Another recent paper that takes into account correlations between the columns of $H$\, is \cite{sel16nonseparable}. A bivariate non-convex penalty is proposed so as to enforce sparsity stronger than alternative convex penalties, while maintaining the convexity of the overall problem. Sparsity within groups is not specifically sought in \cite{sel16nonseparable}.

The penalty proposed in this paper, $P_{\gamma}$ is non-convex. However, its degree of non-convexity is controlled by the parameter $\gamma$ and this in turn allows to formulate convex problems. 
As will be clarified in the sequel (see the proof of Prop.~\ref{prop:weaklyconvex}), $P_{\gamma}$ can be related to the E-Lasso penalty. However, the E-Lasso penalty is convex and can be shown to contain an additive energy term, which in turn penalizes higher coefficients more strongly. Further, the E-Lasso threshold never sets the whole group to zero, unless the group is zero to start with (see \cite{kow09p303}, Remark-6). Thus if a group consists entirely of noise, it will not be totally eliminated, even if it has components with small magnitudes.  The threshold function associated with the proposed penalty function contains a deadzone such that if the coefficients in the group fall in the deadzone, the whole group is eliminated. Therefore the proposed penalty/threshold functions aim to achieve sparsity within and across groups. 

\subsection*{Notation and Preliminaries}
Throughout the paper, vectors are denoted using small case letters, as in $x$. The $i\thh$ component of $x$ is denoted as $x_i$. We are interested in partitions of $x$ into groups in this paper.  We already used $x^{(i)}$ to denote the $i\thh$ subgroup of $x$. That is, for a length-4 vector $x = \begin{pmatrix}x_1, \ldots, x_4 \end{pmatrix}$, if we form two groups of size two, by collecting together consecutive components, we have $x^{(1)} = \begin{pmatrix} x_1, x_2 \end{pmatrix}$\, and $x^{(2)} = \begin{pmatrix} x_3, x_4 \end{pmatrix}$. However, with the exception of Sec.~\ref{sec:hybrid}, the functions under study are separable with respect to groups. Therefore, whenever separability applies, we suppress the group superscript in $x^{(i)}$ and use $x$ to simplify notation, with the understanding that the same discussion applies to all of the groups. 

For a scalar $x \in \mathbb{C}$, the soft threshold function with threshold $\tau > 0$ is defined as,
\begin{equation}
\soft(x,\tau) = \begin{cases}
(|x| - \tau)\,\dfrac{x}{|x|}, & \text{if }|x| > \tau,\\
0, & \text{if }|x| \leq \tau.
\end{cases}
\end{equation}
If $x$\, is a vector, the soft thresholding operator applies to each component of $x$\, separately.

The proximity operator of a convex, lower semi-continuous function $f$\, is defined as \cite{Bauschke,combettes_chp}
\begin{equation}\label{eqn:prox}
J_{\alpha f}(z) = \arg \min_x \frac{1}{2} \| z -x \|_2^2 + \alpha f(x). 
\end{equation}
We also refer to $J$ as the threshold function, if $f$\, under dicsussion is a penalty function.

Throughout the paper, for a given length-$n$ vector $z$ (complex or real valued), we define the cost function $C_{\lambda,\gamma}(x|z)$ as
\begin{equation}\label{eqn:C}
C_{\lambda,\gamma}(x|z) = \frac{1}{2} \| z - x \|_2^2 + \lambda\,P_{\gamma}(x),
\end{equation}
where $P_{\gamma}$\, is given in \eqref{eqn:extendN}.
The threshold function $T_{\lambda,\gamma}(z)$ is defined, in line with \eqref{eqn:prox}, as,
\begin{equation}\label{eqn:T}
T_{\lambda,\gamma}(z) = \arg \min_x C_{\lambda,\gamma}(x|z).
\end{equation}
We  remark that the penalty function used in this paper is not convex but weakly convex. In order for the threshold function to be well-defined, the minimizer of $C_{\lambda,\gamma}(\cdot|z)$ (i.e., the point that minimizes $C_{\lambda,\gamma}(\cdot|z)$) must be unique. To ensure uniqueness, we will check that $C_{\lambda,\gamma}(\cdot|z)$ is strictly convex.

 \subsection*{Outline}
We motivate the proposed penalty function and derive the associated threshold function in Section~\ref{sec:penalty}. We discuss how the non-convex penalty function may be employed to formulate a convex denoising problem with a sparsifying frame and present a minimization algorithm in Section~\ref{sec:denoise}. In Section~\ref{sec:deconvolution}, we present a non-convex deconvolution formulation,  study the convergence of an iterative thresholding algorithm for the presented formulation and demonstrate its performance. Section~\ref{sec:conc} is the conclusion.

\section{A Weakly Convex Penalty}\label{sec:penalty}

The penalty function $\mathbf{P}_{\gamma}$ introduced in \eqref{eqn:Pmain} is separable with respect to the groups and the groups are non-overlapping. Thanks to these properties, it suffices to study the component function $P_{\gamma}(x)$ and the associated threshold function $T_{\lambda,\gamma}$, 
with domain $\mathbb{R}^n$ or $\mathbb{C}^n$. Once  $T_{\lambda,\gamma}$,\, is specified, $\mathbf{T}_{\lambda,\gamma}$\, can be realized by applying $T_{\lambda,\gamma}$ to each group separately.

We start our discussion in Section~\ref{sec:R2}  with penalty/threshold functions defined on $\mathbb{R}^2$,
since this case is easier to visualize and interpret.
After that, we generalize the discussion to $\mathbb{R}^n$ in Section~\ref{sec:Rn}. A discussion of how to tune the parameters and a numerical demonstration of the discussions is provided in Sec.~\ref{sec:parameters}. Extension of the study to $\mathbb{C}^n$ is done in Sec.~\ref{sec:complex}. Finally, in Sec.~\ref{sec:hybrid}, we briefly consider how the proposed penalty can be combined with $\ell_{2,1}$ norms to achieve a modified effect. 

\subsection{The Penalty and the Threshold Function on $\R^2$}\label{sec:R2}

\begin{figure}
\centering
 \includegraphics[scale=1]{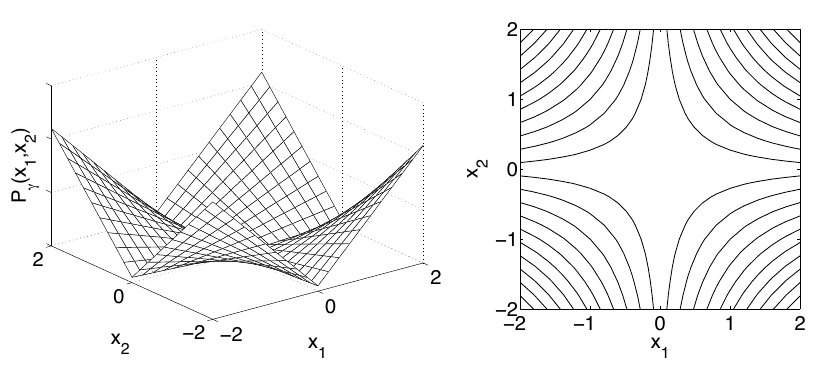}
 \caption{Mesh and contour plots of the proposed $P_{\gamma}(x_1,x_2)$. Notice that the function is not convex. \label{fig:contour}}
\end{figure}

\subsubsection{The Penalty Function}
For a fixed energy vector  $x = (x_1,x_2)$, we seek a penalty function $P$ such that,
\begin{itemize}
\item if $|x_1| \ll |x_2|$ or $|x_2| \ll |x_1|$, $P(x)$ assumes a low value,
\item if $|x_1| \approx |x_2|$, $P(x)$ assumes a high value.
\end{itemize}
Observe that $\bar{P} = |x_1\,x_2|$ satisfies these requirements. However, $\bar{P}$ is exactly zero when one of the components is zero. In order to penalize small non-zero components, we add an $\ell_1$ term and propose the penalty 
\begin{equation}
P_{\gamma}(x) = \gamma\,|x_1\,x_2| + \|x\|_1,
\end{equation}
where $\gamma \geq 0$ is a tuning parameter. Notice that this is the restriction of the function in \eqref{eqn:extendN} to $\mathbb{R}^2$. Mesh and contour plots of this function are shown in Fig.~\ref{fig:contour}. 

 $P_{\gamma}$ is not convex but it becomes convex when we add a quadratic. Such functions are called weakly convex \cite{via83p231}. For $\lambda \leq 1/\gamma$, it can be shown that the function `$\|x\|_2^2/2 + \lambda\,P_{\gamma}(x)$' is convex (see Sec.~\ref{sec:Rn}). Therefore, $\gamma$ may be regarded as a parameter that controls how much $P_{\gamma}$ deviates from being convex. As a consequence of the weak-convexity of $P_{\gamma}$, we find that if $\lambda < 1/\gamma$, for a given $z \in \mathbb{R}^2$, the cost function $C_{\lambda,\gamma}(x|z)$
is strictly convex with respect to $x$.
Thus, $T_{\lambda,\gamma}(z)$ is well-defined when $\lambda\,\gamma <1$. 

Before we further discuss the threshold function, we would like to compare $P_{\gamma}$ to the Elitist-Lasso (E-Lasso) penalty function, which is known to favor large components in a group \cite{kow09p303,kow13p498,kow09p251}. For groups of size two, the E-Lasso penalty uses $P_{\text{EL}}(x)=\|x\|_1^2$. Expanding this, we can write,
\begin{equation}
P_{EL}(x) = 2|x_1\,x_2| + \|x\|_2^2.
\end{equation}
Both $P_{EL}$ and $P_{\gamma}$ employ the term $|x_1\,x_2|$ and this term is responsible for the `elitist' character of the penalties, by which we mean that the penalty assumes lower values for vectors with unequal components. The difference between the two penalties lies in the remaining components. $P_{\gamma}$ contains an additive $\ell_1$ term which helps enforce sparsity if the components have small magnitudes. In contrast, $P_{EL}$ contains an additive energy term, which renders the overall penalty convex at the expense of penalizing high magnitude coefficients more strongly.

\subsubsection{The Threshold Function}
\begin{figure}
\centering
 \includegraphics[scale=1]{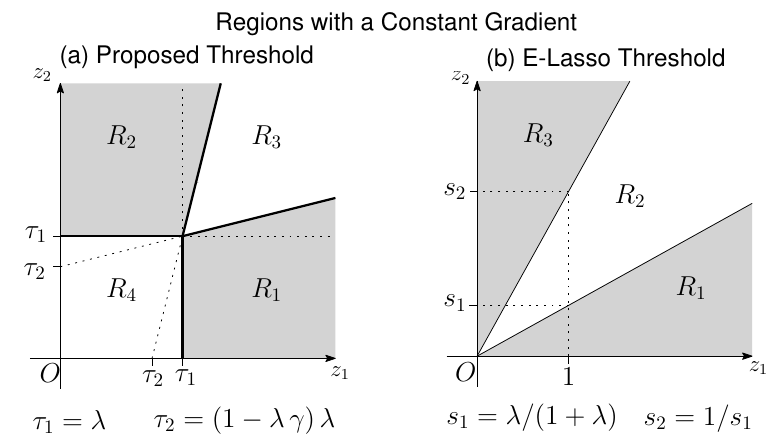}
\caption{The threshold functions associated with the proposed penalty and E-Lasso have constant gradients in the regions indicated above.
\label{fig:regions}}
\end{figure}

The threshold function $T_{\lambda,\gamma}$ can be derived via the optimality conditions for the minimization problem \eqref{eqn:T}.
Let us assume that $z_i \geq 0$, i.e., $z$ lies in the first quadrant (extension to the other quadrants can be achieved by symmetry). Let $\xh = T_{\lambda,\gamma}(z)$. For the regions $R_1$, $R_2$, $R_3$, $R_4$, defined in Fig.~\ref{fig:regions}a, $\xh = (\xh_1,\xh_2)$ is given as follows.
\begin{align*}
&\begin{rcases}
\xh_1 &=  z_1 - \lambda\\
\xh_2 &= 0
\end{rcases},\text{ if } z \in R_1,\\
&\begin{rcases}
\xh_1 &=  0\\
\xh_2 &= z_2 - \lambda
\end{rcases},\text{ if } z \in R_2,\\
&\begin{rcases}
\xh_1  &= \frac{z_1 - \lambda\,\gamma\,z_2 - (1-\lambda\,\gamma)\,\lambda}{1-\lambda^2\,\gamma^2}\\
\xh_2  &= \frac{z_2 - \lambda\,\gamma\,z_1 - (1-\lambda\,\gamma)\,\lambda}{1-\lambda^2\,\gamma^2}
\end{rcases},\text{ if } z \in R_3,\\
&\begin{rcases} \xh_1 &=  0 \\
\xh_2 &= 0
 \end{rcases} \text{ if } z \in R_4.
\end{align*}

\begin{figure}
\centering
 \includegraphics[scale=1]{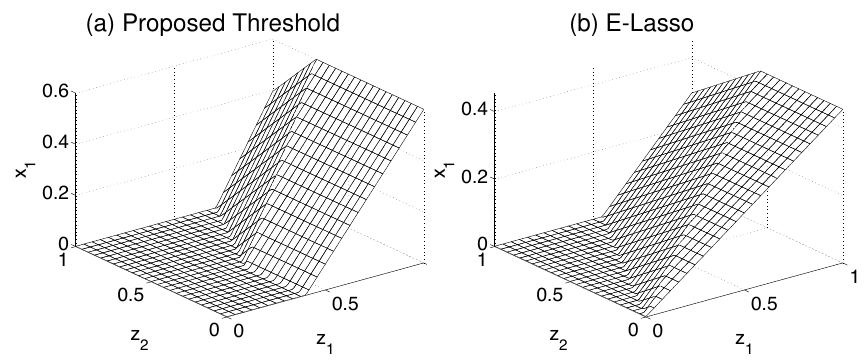}
\caption{The first component of the proposed threshold and E-Lasso.\label{fig:T1}}
\end{figure}

Note that the regions are determined by the two parameters $\lambda$ and $\gamma$ (see $\tau_1$ and $\tau_2$ in Fig.~\ref{fig:regions}a). $R_4$ is the deadzone (i.e., the collection of input vectors which are mapped to zero by the threshold function) for the bivariate threshold function and is determined by the weight $\lambda$. The first component of the threshold function (i.e., the mapping that takes $z = (z_1,z_2)$ to $\xh_1$) is shown in Fig.~\ref{fig:T1}a. Note that for this function, $R_2$ (on which $|z_2| \gg |z_1|$) is also a deadzone.

The proposed threshold function behaves quite differently than a threshold function derived from a separable penalty of the form `$p(x_1) + p(x_2)$'. If the penalty is separable, even if $p$ is non-convex, the deadzone of the  threshold function is rectangular (and hence $\xh_1$ does not depend on $z_2$).

The relevant regions and the first component of the E-Lasso threshold  are shown in Fig.~\ref{fig:regions}b and Fig.~\ref{fig:T1}b respectively. Note that, unlike the proposed threshold (which contains four different regions with constant gradient), the E-Lasso threshold has three regions where its gradient is constant. Also, unlike the proposed threshold, the E-Lasso threshold does not contain a deadzone that eliminates both components.

\subsection{The Penalty and the Threshold Function on $\R^n$}\label{sec:Rn}
We now consider the function $P_{\gamma} : \mathbb{R}^n \to \mathbb{R}_+$ defined in \eqref{eqn:extendN} and the associated threshold function. Unlike $\mathbb{R}^2$, it is not easy to express the threshold function in closed form in $\mathbb{R}^n$. Instead, we will derive a procedure to realize the threshold function. In order to justify the procedure, we will need to first study the properties of the penalty and the threshold functions.

We start with the penalty $P_{\gamma}$. This function is not convex but it becomes convex if we add a quadratic. Such functions are called weakly-convex \cite{via83p231}.
\begin{defn}
For $\gamma \geq 0$, a function $f$ is said to be $\gamma$-weakly convex if 
\begin{equation}
\frac{\gamma}{2}\|x\|_2^2 + f(x)
\end{equation}
is convex.
\end{defn}

\begin{prop}\label{prop:weaklyconvex}
The function $P_{\gamma}$ in \eqref{eqn:extendN} is $\gamma$-weakly convex. Consequently, 
the cost function $C_{\lambda,\gamma}(x|z)$ in \eqref{eqn:C} 
is strictly convex with respect to $x$ if $\lambda\,\gamma < 1$.
\begin{proof}
To see the first claim, observe that,
\begin{equation}
\| x \|_1^2 = \|x\|_2^2 + 2 \sum_{i=1}^{n-1}\,\sum_{m=i+1}^n |x_i\,x_m|.
\end{equation}
Since $\|x \|_1^2$ is convex, this observation implies that the term in the parentheses in \eqref{eqn:extendN} is $1$-weakly convex. Since $\|x\|_1$ is convex, it follows that $P_{\gamma}$ is $\gamma$-weakly convex.

To see the strict convexity of $C_{\lambda,\gamma}$ with respect to $x$, note that it  can be written as the sum of a convex function and 
\begin{equation}\label{eqn:rest}
\frac{1}{2}\|x\|_2^2 + \lambda\,\gamma \, \left( \sum_{i=1}^{n-1}\,\sum_{m=i+1}^n |x_i\,x_m| \right).
\end{equation}
But the function in \eqref{eqn:rest} is strictly convex if $\lambda\,\gamma < 1$ and thus follows the claim.\end{proof}
\end{prop}

It also follows from Prop.~\ref{prop:weaklyconvex} that, when $\lambda\,\gamma < 1$,  $T_{\lambda,\gamma}$ is well-defined thanks to the strict convexity of $C_{\lambda,\gamma}$.

In the following, we will derive two finite-terminating algorithms that realize $T_{\lambda,\gamma}$. For that, we discuss the properties of $T_{\lambda,\gamma}$. We start by showing that $T_{\lambda,\gamma}(z)$ shrinks $z$ towards zero and it is monotone in the sense that it preserves the ordering of $|z_i|$ with respect to $i$. More precisely, we have the following result.
\begin{prop}\label{prop:monotone}
Let $\xh = T_{\lambda,\gamma}(z)$, for $\lambda\,\gamma < 1$.
\begin{enumerate}[(a)]
\item If $z_i \geq 0$, then $z_i \geq \xh_i \geq 0$. If $z_i \leq 0$, then $z_i \leq \xh_i \leq 0$.
\item If $|z_i| > |z_m|$, then $|\xh_i| \geq |\xh_m|$.
\item If $|z_i| = |z_m|$, then $|\xh_i| = |\xh_m|$.
\end{enumerate}
\begin{proof}
See Appendix~\ref{sec:proofmonotone}.
\end{proof}
\end{prop}
We now derive an expression for $T_{\lambda,\gamma}$ using the optimality conditions. Notice that $C_{\lambda,\gamma}$ can be written as,
\begin{multline}
C_{\lambda,\gamma}(x|z) = \frac{1}{2}\|z\|_2^2 -\langle z,x \rangle \\+ \frac{1-\lambda\,\gamma}{2} \|x \|_2^2 + \frac{\lambda\,\gamma}{2} \|x \|_1^2  + \lambda \|x\|_1.
\end{multline}
Using this expression, the optimality conditions are found as,
\begin{equation}
z \in  (1- \lambda\,\gamma) \xh + \lambda\bigl( \gamma \|\xh\|_1 + 1 \bigr) \sign(\xh),
\end{equation}
where $\sign(\xh)$ is a set valued separable function of the vector $\xh$, whose $k\thh$ component is given as,
\begin{equation}\label{eqn:sign}
\sign(\xh)_k = \begin{cases}
\{1\}, &\text{if } \xh_k > 0,\\
[-1,1], &\text{if } \xh_k = 0,\\
\{-1\}, &\text{if } \xh_k < 0.
\end{cases}
\end{equation}

In the following, we assume for simplicity that $z_i$'s are non-negative and ordered, i.e., ${z_1 \geq z_2 \geq \cdots \geq z_n \geq 0}$. The general case can be recovered by a permutation of the vector components and changing signs, thanks to Prop.~\ref{prop:monotone}.  

Prop.~\ref{prop:monotone} implies that there is an index $k \in \{0,1,\ldots, n\}$, such that $\xh_i > 0$ if $i\leq k$ and $\xh_i = 0$ if $i>k$. That is, $k$\, denotes the number of non-zeros in $\xh$. For this special integer $k$, the optimality conditions can be written as,
\begin{subequations}\label{eqn:opt}
\begin{align}
z_i &= (1- \lambda\,\gamma) \xh_i + \lambda\left( \gamma \sum_{m=1}^k \xh_m + 1 \right), \text{if } i \leq k, \label{eqn:nonzero}\\
z_i &\leq \lambda\left( \gamma \sum_{m=1}^k \xh_m + 1 \right), \text{if } i > k.
\end{align}
\end{subequations}
In order to find an expression for $\xh_i$, let us define
\begin{equation}
\bar{z} = \begin{bmatrix}
z_1 \\ z_2 \\ \vdots\\ z_k
\end{bmatrix},
\quad
\bar{x} = \begin{bmatrix}
\xh_1 \\ \xh_2 \\ \vdots\\ \xh_k
\end{bmatrix},
\quad 
\mathbf{1}_k = \begin{bmatrix}
1 \\ 1 \\ \vdots \\1 
\end{bmatrix} \in \mathbb{R}^k.
\end{equation}
We can now express \eqref{eqn:nonzero} as
\begin{equation}\label{eqn:nonzero2}
\bar{z} = (1 - \lambda\,\gamma)\,\bar{x} + \lambda \, \left( \gamma \mathbf{1}^T_k \bar{x} + 1 \right)\,\mathbf{1}_k.
\end{equation}
Multiplying both sides by $\mathbf{1}_k^T$ (noting $\mathbf{1}_k^T\,\mathbf{1}_k =k$) and rearranging, we have
\begin{equation}
\mathbf{1}^T_k \bar{x} = \frac{\mathbf{1}^T_k \bar{z}- k\,\lambda}{1 + (k-1)\lambda\,\gamma}
\end{equation}
Therefore,
\begin{equation}\label{eqn:equality}
\lambda \, \left( \gamma \mathbf{1}^T_k \bar{x} + 1 \right)  = \frac{\lambda\,(1-\lambda\,\gamma) + \lambda\,\gamma\,\sum_{j=1}^k z_j  }{1 + (k-1)\,\lambda\,\gamma}.
\end{equation}
The rhs of \eqref{eqn:equality} will be of interest in the following. 
Let us therefore define, for each $i$,
\begin{equation}\label{eqn:hk}
h(i) = \frac{\lambda\,(1-\lambda\,\gamma) + \lambda\,\gamma\,\sum_{m=1}^i z_m  }{1 + (i-1)\,\lambda\,\gamma}.
\end{equation}
Plugging the expression in \eqref{eqn:equality} back into \eqref{eqn:opt}, we find the equivalent conditions
\begin{subequations}\label{eqn:opt2}
\begin{align}
\xh_i &= (1- \lambda\,\gamma)^{-1} \bigl( z_i - h(k) \bigr), &\text{if } i \leq k,\\
z_i &\leq h(k), &\text{if } i > k.
\end{align}
\end{subequations}
Notice that the requirement $\xh_i > 0$ for $i\leq k$ implies that $z_i > h(k)$ for $i\leq k$.
The foregoing discussion is summarized in the following proposition. 
\begin{prop}\label{prop:thold}
Let $\xh = T_{\lambda,\gamma}(z)$, for $\lambda\,\gamma < 1$. Also, let $k$ denote the number of non-zeros of $\hat{x}$. Then,
\begin{equation}\label{eqn:xtz}
\xh = (1 - \lambda\,\gamma)^{-1}\,\soft\bigl(z, h(k)\bigr),
\end{equation}
where 
\begin{equation}
h(k) = \frac{\lambda\,(1-\lambda\,\gamma) + \lambda\,\gamma\,\sum_{m=1}^k |z_m|  }{1 + (k-1)\,\lambda\,\gamma},
\end{equation}
(with the convention $\sum_{m=1}^0 |z_m| = 0$).
\end{prop}
We remark that the description of the threshold function in Prop.~\ref{prop:thold} is implicit because the integer $k$\, in \eqref{eqn:xtz}, namely the number of non-zeros in $\hat{x}$, depends on $\hat{x}$. 
We next discuss how to determine the integer $k$.
We will present two different search schemes for finding the correct value of $k$. The following lemma will be useful for that end.
\begin{lem}\label{lem:selectk}
Suppose $z_1 \geq z_2 \geq \cdots z_n \geq 0$ and $\lambda\,\gamma < 1$. Let $h(i)$ be defined as in \eqref{eqn:hk}. Then, 
\begin{enumerate}[(a)]
\item if $z_{i+1} > h(i)$,  then $z_i > h(i-1)$,
\item if $z_{i} \leq h(i)$,  then $z_{i+1} \leq h(i+1)$.
\end{enumerate}
\begin{proof}
See Appendix~\ref{sec:proofselectk}.
\end{proof}
\end{lem}

We can use this lemma to develop a procedure for determining $k$. It follows from the lemma that we can start from $k= 0$ and keep increasing $k$ until $h(k) > z_{k+1}$. This procedure is summarized in Algorithm~\ref{algo:klinear}.
\begin{algorithm}[h!]{\caption{Realization of $T_{\lambda,\gamma}$ -- Linear Search for $k$}\label{algo:klinear}}
\begin{algorithmic}[0]
\REQUIRE $y\in \mathbb{R}^n$
\STATE $z \gets \text{descending-sort}(|y|)$
\STATE $k\gets 0$
\WHILE{$h(k) < z_{k+1}$}
\STATE $k \gets k+1$ \COMMENT{increment $k$}
\ENDWHILE
\STATE $\xh\gets (1-\lambda\,\gamma)^{-1}\,\soft\bigl(y, h(k) \bigr)$.
\end{algorithmic}
\end{algorithm}

If $k$\, is suspected to be small, then this algorithm terminates quickly. In the worst case, the algorithm will execute the `while' loop $n$ times. If, however, $n$ is large and $k$ is not expected to be small, then a binary search for $k$ might be computationally more suitable. The following discussion, that relies on Lemma~\ref{lem:selectk} implies that such a binary search terminates.

Suppose we pick an arbitrary $i$ and check the following conditions.
\begin{subequations}
\begin{align}
z_i &> h(i),\label{eqn:cond1} \\
z_{i+1} & \leq h(i). \label{eqn:cond2}
\end{align}
\end{subequations}
Notice that since $z_i \geq z_{i+1}$, the conditions cannot be violated simultaneously. Now observe that
\begin{itemize}
\item If both conditions hold, the current guess of $i$ is equal to the sought $k$.
\item If \eqref{eqn:cond1} holds and \eqref{eqn:cond2} is violated, then by Lemma~\ref{lem:selectk},  $k$ must be greater than $i$.
\item If \eqref{eqn:cond2} holds and \eqref{eqn:cond1} is violated, then again by Lemma~\ref{lem:selectk}, $k$ must be less than $i$.
\end{itemize} 

These observations lead to an implementation of $T_{\lambda,\gamma}$ as in Algorithm~\ref{algo:nThold}. In contrast to the $O(n)$ complexity of Algorithm~\ref{algo:klinear}, this algorithm has $O\bigl(\log(n)\bigr)$ complexity.

\begin{algorithm}[h!]{\caption{Realization of $T_{\lambda,\gamma}$ -- Binary Search for $k$}\label{algo:nThold}}
\begin{algorithmic}[0]
\REQUIRE $y\in \mathbb{R}^n$
\STATE $z \gets \text{descending-sort}(|y|)$
\STATE $\text{flag} \gets \text{true}$
\IF{$z_1 < \lambda$}
\STATE $k\gets 0$ 
\STATE $\text{flag} \gets \text{false}$
\ELSIF{$z_n > h(n)$}
\STATE $k\gets n$
\STATE $\text{flag} \gets \text{false}$
\ELSE
\STATE $k_0 \gets 0$  \COMMENT{left end of the interval}
\STATE $k_1 \gets n$  \COMMENT{right end of the interval}
\ENDIF
\WHILE{flag}
\STATE $k \gets \lfloor (k_0 + k_1)/2 \rfloor$ \COMMENT{middle of the working interval}
\IF{$z_k > h(k)$ and $z_{k+1} \leq h(k)$}
\STATE $\text{flag} \gets \text{false}$ \COMMENT{correct value of $k$ is found}
\ELSIF{$z_k \leq h(k)$}
\STATE $k_1 \gets k$ \COMMENT{update the right end}
\ELSE
\STATE $k_0 \gets k$ \COMMENT{update the left end}
\ENDIF
\ENDWHILE
\STATE $x \gets (1-\lambda\,\gamma)^{-1}\,\soft\bigl(y, h(k)\bigr)$.
\end{algorithmic}
\end{algorithm}

\subsection{Tuning the Parameters of the Threshold Function}\label{sec:parameters}

Let us now discuss some special cases to better understand the role of the parameters $\lambda$ and $\gamma$. As in the previous subsection, we will assume that $z_1 \geq \cdots \geq z_n \geq 0$ and $\xh = T_{\lambda,\gamma}(z)$.

Observe first that $h(0) = \lambda$. If $z_i< \lambda$ for all $i$, then $\xh = 0$. Thus the deadzone of the threshold function is a cube of width $\lambda$ in $\mathbb{R}^n$.

Suppose now that $z_1 > \lambda$. In that case, we will definitely have $\xh_1 >0$. We find that
\begin{equation}
h(1) = \lambda + \lambda\gamma (z_1 - \lambda).
\end{equation}
Notice that in order for $\xh_2$ to be non-zero, the threshold that $z_2$ needs to exceed has increased from $\lambda$\, by an amount proportional to $(z_1 - \lambda)$. The higher $z_1$ is, the higher will be the new threshold. In fact, observe that as $\lambda\,\gamma \to 1$, the threshold converges to $z_1$. Since $z_2 \leq z_1$, we can therefore force only a single component to survive by choosing $\gamma$ close to $1/\lambda$. When $z_2 < h(1)$, we find that,
\begin{equation}
\xh_1 = z_1 - \lambda.
\end{equation}
Thus the single surviving component is obtained by soft thresholding the largest component with $\lambda$. 

The following proposition provides further information on how the potential thresholds $h(i)$\, behave for arbitrary $i$.
\begin{prop}\label{prop:hi}
Suppose $z_1 \geq z_2 \geq \cdots z_n \geq 0$ and $\lambda\,\gamma < 1$. Let $h(i)$ be defined as in \eqref{eqn:hk}. If $z_{i+1} > h(i)$, then ${z_{i+1} > h(i+1) > h(i)}$.
\begin{proof}
See Appendix~\ref{sec:proofhi}.
\end{proof}
\end{prop}

We know that if $z_{i+1} > h(i)$, then $h(i)$ is not the actual threshold and  $k>i$. Prop.~\ref{prop:hi} implies that $h(k)$\, is actually greater than $h(i)$, but it is bounded above by $z_{i+1}$. In fact, we can deduce from Prop.~\ref{prop:hi} that 
\begin{equation}
z_1  \geq \cdots \geq z_k > h(k) \geq h(k-1 ) \geq \cdots \geq h(0) = \lambda.
\end{equation}
In the case where the observations are purely noise, we would like to set $\xh_i=0$ for all $i$. This motivates choosing $\lambda = c\sigma$, where $\sigma$ denotes the noise standard deviation and $c$ is a constant around unity. Once we fix the value of $\lambda$, the number of non-zero components, $k$, and the threshold $h(k)$ will depend on $\gamma$ (and $z$). The following proposition provides precise bounds on $\gamma$.

\begin{prop}\label{prop:gamma}
Let $\xh = T_{\lambda,\gamma}(z)$ where $\lambda\,\gamma < 1$ and $z_i \geq 0$ for all $i$.
${\xh_1 \geq  \cdots \geq \xh_k >0}$ and ${\xh_{k+1}  = \cdots = \xh_n = 0}$\, if and only if 
\begin{subequations}\label{eqn:gamma}
\begin{align}
\lambda\,\gamma &> \frac{(z_{k+1} - \lambda )_+}{(z_{k+1} - \lambda)_+ + \sum_{i=1}^k\,(z_i - z_{k+1})},\label{eqn:5a}\\
\lambda\,\gamma &< \frac{(z_k - \lambda )_+}{(z_k - \lambda)_+ + \sum_{i=1}^{k-1}\,(z_i - z_k)}\label{eqn:5b}.
\end{align}
\end{subequations}
\begin{proof}
See Appendix~\ref{sec:proofgamma}.
\end{proof}
\end{prop}

\begin{figure}
\centering 
\includegraphics[scale=1]{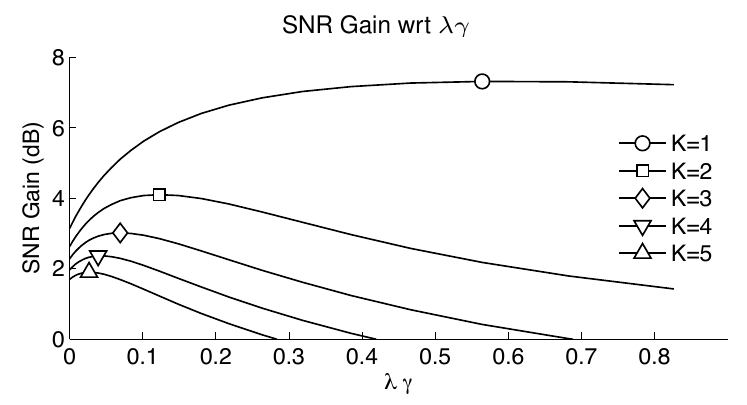}
\caption{Average SNR gains with respect to $\lambda,\gamma$, where $\lambda$ is fixed and $\gamma$ is varied. $K$ denotes the number of non-zero components present in the clean signal.\label{fig:gamma}}
\end{figure}
Prop.~\ref{prop:gamma} suggests that, if we would like to retain more components in the estimate $\xh$, then we need to choose a small $\gamma$ so that $\lambda \gamma$ is small. Also, if the signal is scaled by multiplying with a factor $\beta$, then the numbers on the rhs in \eqref{eqn:gamma} stay approximately the same. Therefore, the constant $\lambda\,\gamma$ controls the  number of non-zeros, almost independently of the scale. To summarize, the following may serve as a guide for the selection of the parameters $\lambda$ and $\gamma$. 
\begin{itemize}
\item $\lambda$\, can be chosen  proportional to the noise standard deviation.
\item Once $\lambda$ is chosen, $\gamma$ can be selected independently of the scale of the signal but should satisfy $\lambda\,\gamma < 1$. The product $\lambda\,\gamma$ determines the number of non-zeros in the estimate. Higher it is, lower the number of non-zeros will be.
\end{itemize}

In order to demonstrate the foregoing discussion, we now consider two simple experiments on synthetic signals.

In both experiments, the desired signal is of length 10 and it has $K$ non-zero components, where the non-zero values are obtained by sampling from a Gaussian distribution. We add Gaussian noise to this signal so that the observation SNR is 5~dB.

For the first experiment, we apply $T_{\lambda,\gamma}$ to the observation for a fixed $\lambda = \sigma /2$ and varying $\gamma$. We repeat the experiment for 10K trials to obtain an average figure. The average gain in SNR (dB) with respect to $\lambda\,\gamma$ is shown in Fig.~\ref{fig:gamma}. We see from the figure that the best $\lambda\,\gamma$ value decreases with increasing $K$. This is in line with Prop.~\ref{prop:gamma}, which suggests that in order for the reconstruction to have more non-zeros, the product $\lambda\,\gamma$ must be smaller.

\begin{figure}
\centering
 \includegraphics[scale=1]{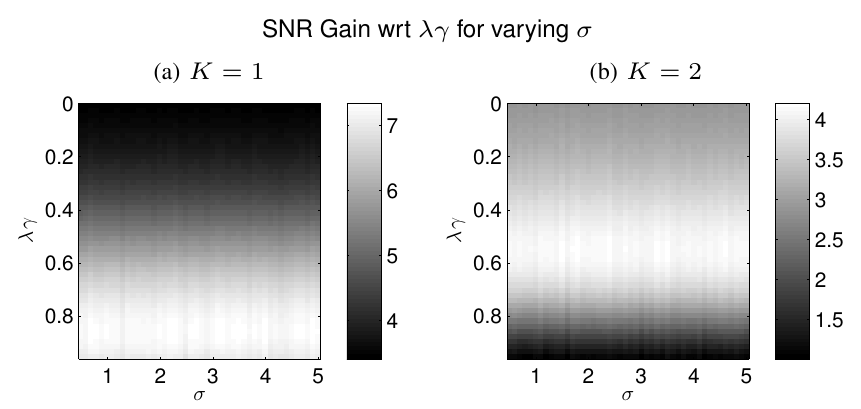}
\caption{Average SNR gains with respect to the product $\lambda\gamma$, where $\lambda$ is fixed and $\gamma$ is varied. $K$ denotes the number of non-zero components present in the clean signal. As $\sigma$\, varies, the value of $\lambda \gamma$\, that maximizes the SNR stays roughly constant.\label{fig:ConstGam}}
\end{figure}

In a second experiment, we demonstrate that the best choice of the product $\lambda \gamma$\,  is not affected much by the scale of the observations. For this, we vary the standard deviation of the noise, $\sigma$,\, as well as the signal energy so that the SNR stays constant at 5~dB. We set $\lambda = \sigma /2$ as in the previous experiment. For fixed $K$ value and varying $\sigma$,  we apply the threshold using a range of $\gamma$\, values satisfying $\lambda\,\sigma < 1$ and compute the gain in SNR. For each value of $\sigma$\, and $\gamma$, an average SNR gain value is obtained by repeating this for 10K trials. The resulting SNR gain images are shown in Fig.~\ref{fig:ConstGam}, for $K=1$ and $K=2$. Observe that for both $K=1$\, and $K=2$, the best choice of $\lambda\,\gamma$\, stays approximately constant as $\sigma$\, varies. 
\subsection{Extension to $\mathbb{C}^n$}\label{sec:complex}

For $x\in \mathbb{C}^n$, we extend $P_{\gamma}$ straightforwardly as,
\begin{equation}
P_{\gamma}^c(x) = \gamma \, \left( \sum_{i=1}^{n-1}\,\sum_{j=i+1}^n |x_i\,x_j| \right) + \|x\|_1.
\end{equation}
The threshold function is similarly defined as,
\begin{equation}\label{eqn:TholdComplex}
T_{\lambda,\gamma}^c(z) = \arg \min_{x \in \mathbb{C}^n} \, \frac{1}{2} \|x - z \|_2^2 + \lambda\,P^c_{\gamma}(x).
\end{equation}
Fortunately, the threshold function derived for $\mathbb{R}^n$ applies for $\mathbb{C}^n$ with a little modification. The following observation is useful for showing that.
\begin{lem}\label{lem:arg}
Suppose $z\in \mathbb{C}^n$ and $\xh = T_{\lambda,\gamma}^c(z)$. If $|\xh_i|>0$, then $\arg(\xh_i) = \arg(z_i)$.
\begin{proof}
Suppose $\arg(\xh_i) \neq \arg(z_i)$. Define $\tilde{x}$  such that $|\tilde{x}_i| = |\xh_i|$ for all $i$ and for $|\tilde{x}_i|>0$, set $\arg(\tilde{x}_i) = \arg(z_i)$. Then, $P_{\gamma}^c(\tilde{x}) = P_{\gamma}^c(\xh)$ but $\|z - \tilde{x}\|_2^2 < \| z- \xh \|_2^2$, contradicting the fact that $\xh$ minimizes the cost in \eqref{eqn:TholdComplex}.
\end{proof}
\end{lem}
With the help of this lemma, we obtain an expression for $T_{\lambda,\gamma}^c$ in terms of $T_{\lambda,\gamma}$.
\begin{prop}
Suppose $z\in \mathbb{C}^n$. Let $|z|$ denote the vector containing the magnitudes of the components of $z$. Let ${\xh = T_{\lambda,\gamma}^c(z)}$ and $u = T_{\lambda,\gamma}(|z|)$. Then, $\xh_k  = u_k\,e^{j\arg(z_k)}$.
\begin{proof}
Notice that $|z_k| = z_k\,e^{-j\arg(z_k)}$. Using this observation, it can be shown by a change of variables that if $\tilde{x} = T_{\lambda,\gamma}^c(|z|)$, then  $\tilde{x}_k = x_k\,e^{-j\arg(z_k)}$. Now since $\arg(|z_k|) = 0$, for all $k$, it follows by Lemma~\ref{lem:arg} that $\tilde{x}_k$ are real and non-negative. Thus, for the input $|z|$, we can restrict the minimization in \eqref{eqn:TholdComplex} to $\mathbb{R}^n$. Thus $\tilde{x} = T_{\lambda,\gamma}(|z|) = u$ and the claim follows.
\end{proof}
\end{prop}
It follows from this proposition that the threshold function on $\mathbb{C}^n$ can be realized by applying $T_{\lambda,\gamma}$ to the magnitudes of the input, followed by a correction of the argument of the complex number. For this reason, in the following, we will not differentiate between $T_{\lambda,\gamma}$  and $T_{\lambda,\gamma}^c$. 

\subsection{An Extension to Sub-Groups}\label{sec:hybrid}
We have so far considered a vector $x \in \mathbb{C}^n$ to comprise a group belonging to a larger signal. We now add an additional layer and consider subgroups of $x$\, to define a hybrid penalty, that can be used to complement $\ell_{2,1}$\, norms. In this setting, we refer to $x$ as a `super-group'. Specifically, suppose $x$\, is partitioned into $m$\, non-overlapping  sub-groups, i.e. $x = \begin{bmatrix} x^{(1)} & x^{(2)} & \ldots & x^{(m)} \end{bmatrix}$. Also, let $w$\, denote the length-$m$ vector whose $i\thh$ component is $w_i = \|x^{(i)}\|_2$. We define a hybrid penalty $\tilde{P}_{\gamma}(x)$ as,
\begin{equation}
\tilde{P}_{\gamma}(x) = P_{\gamma}(w).
\end{equation}
Observe that, 
\begin{align}\label{equation}
\tilde{P}_{\gamma}(x) &= \frac{\gamma}{2}\bigl( \| w\|_1^2 -  \|w\|_2^2 \bigr) + \| w \|_1 \\
& = \frac{\gamma}{2} \, \bigl( \| x\|_{2,1}^2 -  \|x\|_2^2 \bigr) + \|x\|_{2,1},
\end{align}
where $\|x\|_{2,1} = \sum_{i=1}^m \| x^{(i)}\|_2$. Therefore, $\tilde{P}_{\gamma}$ is $\gamma$-weakly convex. The threshold function of $\tilde{P}_{\gamma}$ is similarly defined as 
\begin{equation}\label{eqn:Ttilde}
\tilde{T}_{\lambda,\gamma}(z)  = \arg \min_{x} \frac{1}{2} \| x - z \|_2^2 + \lambda\,\tilde{P}_{\gamma}(x).
\end{equation}
Thanks to the weak-convexity of $\tilde{P}_{\lambda,\gamma}$, $\tilde{T}_{\lambda,\gamma}$ is well-defined for $\lambda\,\gamma <1$.
$\tilde{T}_{\lambda,\gamma}$ can be easily described using $T_{\lambda,\gamma}$, as follows (see also \cite{kow13p498} for a discussion of convex multi-layer penalties).
\begin{prop}\label{prop:hybrid}
Suppose a vector $z$ partitioned into $m$ groups where the $i\thh$\, group is denoted as $z^{(i)}$. Also, let $w$ denote the length-$m$ vector whose $i\thh$\, component is $w_i = \| z^{(i)}\|_2$. Let $\hat{x} = \tilde{T}_{\lambda,\gamma}(z)$, $\hat{w} = T_{\lambda,\gamma}(w)$. If we partition $\hat{x}$ into $m$ groups similarly as $z$ (with $\hat{x}^{(i)}$ denoting the $i\thh$ group),  we have,
\begin{equation}
\hat{x}^{(i)} = \begin{cases}
\dfrac{\hat{w}_i}{w_i}\, z^{(i)}, &\text{ if } w_i > 0,\\
0, &\text{ if } w_i = 0.
\end{cases}
\end{equation}
\begin{proof}
See Appendix~\ref{sec:proofhybrid}.
\end{proof}
\end{prop}

In words, this proposition implies that the orientation of $x^{(i)}$ is the same as that of $z^{(i)}$. To find the length of $\xh^{(i)}$, i.e., $\|\xh^{(i)}\|_2$,  we apply the thresholding operator $T_{\lambda,\gamma}$ to $w$.
We will consider an application of this penalty and threshold function in Sec.~\ref{sec:denoise} for a convex denoising formulation.

\section{Application-I : Convex Denoising with a Sparsifying Frame}\label{sec:denoise}

We now consider the application of the proposed penalty in a denoising problem, when a sparsifying frame is given. That is, we assume that we have available a linear transform with a stable inverse (see \cite{Christensen} for a detailed discussion) which allows to represent the signal with high fidelity using a small number of transform domain coefficients. We will specifically seek a convex formulation for this problem. 

\subsection{A Convex Denoising Formulation}\label{sec:AnalysisSynthesis}

Let $y$ be a noisy observation of a clean signal $x_c$ for which a sparsifying frame is given. Let $S$ and $S^*$ denote the analysis and synthesis operators for the frame \cite{Christensen}. We assume that $S^*\,S = I$, i.e., the frame is Parseval \cite{Christensen}. We have two choices for formulating the denoising problem, namely synthesis and analysis prior formulations \cite{ela07p947,selFig09SPIE}. The two behave quite differently under a non-convex penalty such as the one considered in this paper.

In the setting described above, the synthesis prior denoising formulation is,
\begin{equation}\label{eqn:synth}
\min_t \, \frac{1}{2}\,\|y - S^*\,t\|_2^2 + \lambda\,\mathbf{P}_{\gamma}(t).
\end{equation}
If we denote the minimizer as $\hat{t}$, the denoised estimate is given as $\hat{x} = S^*\,\hat{t}$.
In order to investigate convexity, let us rewrite the cost function in \eqref{eqn:synth} as 
\begin{equation}\label{eqn:synth2}
\left[ \frac{1}{2}\,\|y - S^*\,t\|_2^2  - \frac{\alpha}{2} \| t\|_2^2 \right] + \left[ \frac{\alpha}{2}\|t\|_2^2 + \lambda\,\mathbf{P}_{\gamma}(t) \right].
\end{equation}
Notice that the second term in square brackets in \eqref{eqn:synth2} is convex if $\lambda\,\gamma \leq \alpha$. The Hessian of the first term in square brackets is $S\,S^* - \alpha I$. Thus, the first term is also convex if $S\,S^* \succeq \alpha I$. In that case, the problem in \eqref{eqn:synth} will be convex. However, if the frame is overcomplete, the condition $S\,S^* \succeq \alpha I$ is not satisfied for $\alpha >0$. 
Therefore, we can guarantee the convexity of the synthesis prior problem only for $\gamma = 0$, for which $\mathbf{P}_{\gamma}$ is equivalent to an $\ell_1$ norm. This leads us to consider the analysis prior formulation given as,
\begin{equation}\label{eqn:analysis}
\min_x \, \frac{1}{2}\,\|y - x\|_2^2 + \lambda\,\mathbf{P}_{\gamma}(S\,x).
\end{equation}
\begin{prop}
Suppose $S$ is the analysis operator of a Parseval frame. The problem in \eqref{eqn:analysis} is convex if $\lambda\,\gamma \leq 1$.
\begin{proof}
Since the frame is tight, we have $\|Sx\|_2^2 = \|x\|_2^2$. Therefore the cost function in \eqref{eqn:analysis} can be written as,
\begin{multline}\label{eqn:f1f2}
\frac{1}{2}\,\|y - x\|_2^2 - \frac{1}{2}\|x\|_2^2 + \frac{1}{2}\|Sx\|_2^2 + \lambda\,\mathbf{P}_{\gamma}(S\,x) 
\\
=\underbrace{ \frac{1}{2}\,\|y\|_2^2 - \langle x,y\rangle}_{f_1(x)} + \underbrace{ \frac{1}{2}\|Sx\|_2^2 + \lambda\,\mathbf{P}_{\gamma}(S\,x)}_{f_2(Sx)}.
\end{multline}
In \eqref{eqn:f1f2}, $f_1$ is an affine function and is therefore convex. The function $f_2(x)$ in \eqref{eqn:f1f2} is convex when $\lambda\,\gamma \leq 1$, by Prop.~\ref{prop:weaklyconvex}. Since pre-composition with a linear operator preserves convexity \cite{HiriartFund}, $\tilde{f}(x) = f_2(Sx)$ is also convex. Thus the cost function in \eqref{eqn:analysis} can be expressed as the sum of two convex functions and is therefore convex.
\end{proof}
\end{prop}

\subsubsection{The Douglas-Rachford Algorithm}
In order to obtain a minimizer of \eqref{eqn:analysis}, we adapt the Douglas-Rachford algorithm \cite{lio79p964,eck92p293,combettes_chp}.
The Douglas-Rachford algorithm is suitable for minimization problems of the form
\begin{equation}\label{eqn:DRproblem}
\min_t f(t) + g(t),
\end{equation}
where both $f$ and $g$ are convex. The Douglas-Rachford iterations for such a problem are,
\begin{equation}\label{eqn:DR}
t^{k+1} = \frac{1}{2} t^k + \frac{1}{2}\bigl( 2J_{\alpha f} - I \bigr) \bigl( 2J_{\alpha g} - I \bigr)\,t^k,
\end{equation}
where $\alpha >0$ is a parameter and $J_{\alpha f}$, $J_{\alpha g}$ are proximity operators associated with $f$ and $g$ (as defined in \eqref{eqn:prox}). 
The sequence $t^k$ constructed in \eqref{eqn:DR} converges to some $t^*$ such that $J_{\alpha g}(t^*)$ minimizes \eqref{eqn:DRproblem}.
\subsubsection{Adapting the Douglas-Rachford Algorithm}
The problem in \eqref{eqn:analysis} is not readily suitable for the application of the Douglas-Rachford algorithm. We now transform the problem to write it in a suitable form.

Since $S^*\,S = I$, we have  
\begin{equation}
\|x -y \|_2^2 = \| Sy - Sx\|_2^2.
\end{equation}
Now if $\mathcal{R}(S)$ denotes the range of $S$, we can change variables and obtain a problem equivalent to \eqref{eqn:analysis} 
 as,
\begin{equation}\label{eqn:equivalent}
\min_u \, \underbrace{ \frac{1}{2}\,\|S\,y - u \|_2^2 + \lambda\,P_{\gamma}(u)}_{f(u)} + \underbrace{i_{\mathcal{R}(S)}(u)}_{g(u)},
\end{equation}
where $i_{\mathcal{R}(S)}$ is the indicator function of $\mathcal{R}(S)$ \cite{HiriartFund}. If  $u^*$ denotes a minimizer of \eqref{eqn:equivalent}, then $S^*\,u^*$ minimizes  \eqref{eqn:analysis}. 
In this formulation, both $f$ and $g$ are convex, provided that $\lambda\,\gamma \leq 1$. Thus the Douglas-Rachford algorithm is applicable for this splitting. We remark that in this setting, the proximity operator for $g = i_{\mathcal{R}(S)}$ is simply a projection onto $\mathcal{R}(S)$ (see e.g. \cite{combettes_chp}), which can be achieved by applying $S\,S^*$, thanks to the Parseval property of the frame. The proximity operator for $f$ can be expressed in terms of the threshold function as follows.
\begin{subequations}
\begin{align}
&J_{\alpha f}(z) \nonumber \\
&= \arg \min_u \frac{1}{2\alpha} \| z - u\|_2^2 + \frac{1}{2}\|Sy - u\|_2^2  + \lambda\,\mathbf{P}_{\gamma}(u) \label{eqn:a}\\
&=\arg \min_u \frac{1}{2}\left\| u - \frac{\alpha}{\alpha+1}\left(Sy + \frac{z}{\alpha} \right) \right\|_2^2 +  \frac{\alpha}{\alpha+1}\lambda\,\mathbf{P}_{\gamma}(u) \label{eqn:b}\\
&= \mathbf{T}_{(\beta\,\lambda), \gamma}\biggl( \beta\,\left(Sy + \frac{z}{\alpha} \right) \biggr),
\end{align}
\end{subequations}
where $\beta = \dfrac{\alpha}{\alpha+1}$. We remark that in passing from \eqref{eqn:a} to \eqref{eqn:b}, we discarded an additive term and removed a positive factor from the cost function (equalities remain valid because we are seeking the minimizer).

Resulting pseudocode for the Douglas-Rachford iterations for this problem is given in Algorithm~\ref{algo:Denoise}.

\begin{algorithm}[h!]{\caption{Analysis Prior Denoising Algorithm}\label{algo:Denoise}}
\begin{algorithmic}
\STATE Initialize $\alpha >0$, $t$.
\STATE Set $\beta \gets (1+\alpha)^{-1}\,\alpha $
\REPEAT
\STATE $u \gets S\,S^*\,t$
\STATE $z  \gets \mathbf{T}_{\beta\,\lambda,\gamma}\Bigl(\beta \bigl(Sy + \alpha^{-1}(2u - t)\bigr)\Bigr)$
\STATE $t  \gets z + t - u$
\UNTIL{convergence}
\STATE $x^* \gets S^*\,t$
\end{algorithmic}
\end{algorithm}

\subsection{Numerical Experiment}\label{sec:DenoiseExperiment}
\begin{figure}
\centering
 \includegraphics[scale=1]{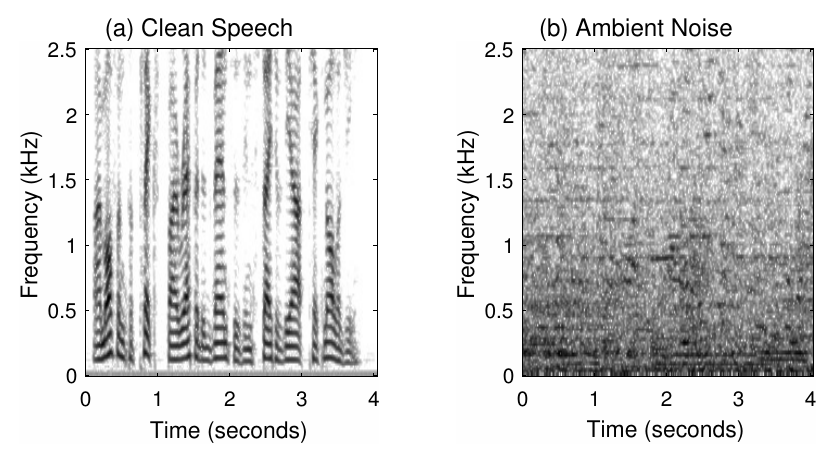}
\caption{Spectrograms of (a) the clean speech and  (b) the noise signal used in the analysis prior denoising experiment.\label{fig:Exp1in}}
\end{figure}

We now demonstrate how the denoising formulation/algorithm performs on an audio signal and compare it against formulations that employ different regularizers. The clean signal is a speech signal, sampled at 16~kHz, whose spectrogram is shown in Fig.~\ref{fig:Exp1in}a. We use the short-time Fourier transform (STFT) as the tight frame in this experiment. The window size is 60~ms (960 samples) and the hop-size is 15~ms (240 samples). 

For regularization, we compare the $\ell_1$\, norm, E-Lasso, the $\ell_{2,1}$ norm and two different versions of the proposed penalty. Our aim in this experiment is two-fold. First, we demonstrate the difference between the proposed penalty, E-Lasso and the $\ell_1$ norm. Second, we show that the proposed penalty can be used to complement the $\ell_{2,1}$ norm to obtain enhanced reconstructions. 

\begin{figure}
\centering
 \includegraphics[scale=1]{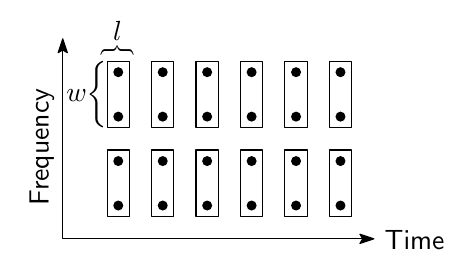}
\caption{Definition of the groups in the time-frequency lattice used in the denoising experiment in Section~\ref{sec:DenoiseExperiment}.
The parameters $l$ and $w$ denote the length along the time axis and the width along the frequency axis respectively. 
}\label{fig:TFGroups}
\end{figure}

In order to describe the group penalties, consider Fig.~\ref{fig:TFGroups}, which introduces notation for $P_{EL}$ (E-Lasso penalty), $P_\gamma$ (proposed)  and $\ell_{2,1}$ norm. For $P_{EL}$ and $P_{\gamma}$, we select the length along the time axis as $l=1$ and the width along the frequency axis as $w = 16$. This covers a frequency bandwidth of 320~Hz. Our aim is to exploit the isolated appearance of the harmonics viewed along the frequency axis. In contrast, for the $\ell_{2,1}$ norm, we take $w = 1$, $l=8$. With this choice, we aim to collect the coefficients belonging to a harmonic into a single group. However, unlike $P_{EL}$ and $P_{\gamma}$, $\ell_{2,1}$ norm regularization does not specifically seek to isolate the harmonics like  $P_{\text{EL}}$ and $P_{\gamma}$. In order to obtain such an effect, we add an additional layer of grouping as depicted in Fig.~\ref{fig:TFGroups} and use the penalty $\tilde{P}_{\gamma}$ introduced in Sec.\ref{sec:hybrid}. We stack 16 neighboring groups along the frequency axis, used in the $\ell_{2,1}$ norm to define non-overlapping super-groups for $\tilde{P}_{\gamma}$.

\begin{figure}
\centering
 \includegraphics[scale=1]{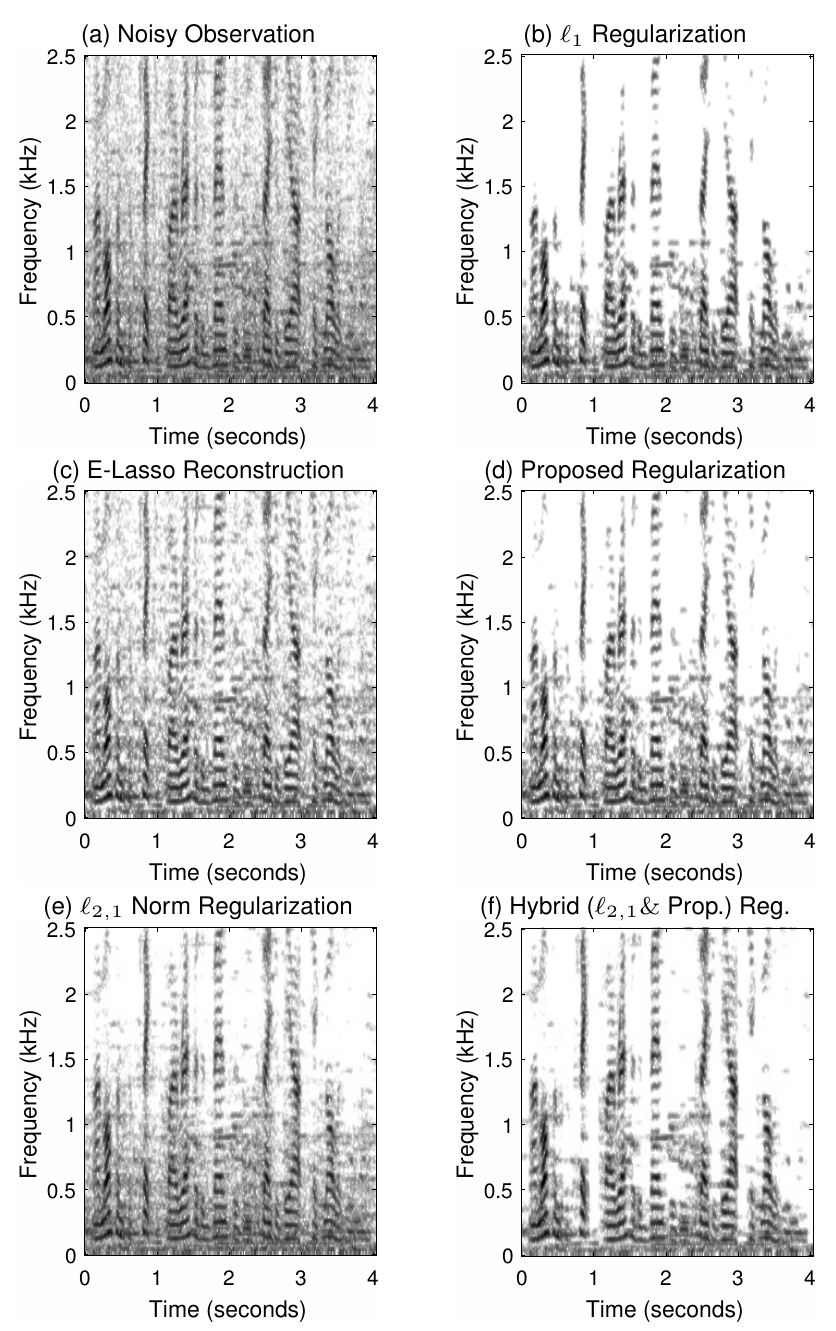}
\caption{Spectrograms of (a) the noisy signal  and the reconstructions using (b)~the $\ell_1$ norm, (c)~the E-Lasso penalty, (d)~the proposed penalty, (e)  the $\ell_{2,1}$ norm, (f) the hybrid penalty in Sec.~\ref{sec:hybrid} obtained by combining the $\ell_{2,1}$ norm with the proposed penalty. \label{fig:Exp1Out}}
\end{figure}

We produce a noisy observation by adding noise (see Fig.~\ref{fig:Exp1in}b) consisting of ambient sounds recorded in a casino (machine sounds and crowd noise). Notice that energy of the ambient noise is not uniform over the frequencies and decreases with increasing frequency. Therefore, this noise can be considered pink. The input SNR is 5~dB. The spectrogram of the noisy signal is shown in Fig.~\ref{fig:Exp1Out}a. For this observation, we perform denoising using the different regularizers in the analysis prior formulation, namely \eqref{eqn:analysis}. The denoising algorithms for the different regularizers can be obtained by replacing $\mathbf{T}_{\lambda,\gamma}$ with the corresponding threshold functions in Algorithm~\ref{algo:Denoise}. The value of the regularizer weight $\lambda$ is selected by a sweep search for the $\ell_1$ norm,  E-Lasso and the $\ell_{2,1}$ norm (output SNR maximizing value is chosen). For $P_{\lambda,\gamma}$, we set $\lambda$ to be half the value of $\lambda$ used for the $\ell_1$ norm and perform a sweep search for $\gamma$ subject to $\gamma < 1/\lambda$, to maximize the output SNR. For $\tilde{P}_{\gamma}$, we similarly set $\lambda$\, equal to half the value used for the $\ell_{2,1}$ norm and search for the best $\gamma$. The optimal choices of $\gamma$\, were found as 2.68 and 7.85 for $P_{\gamma}$ and $\tilde{P}_{\gamma}$, respectively.

The resulting reconstructions are shown in Fig.~\ref{fig:Exp1Out}. The output SNRs  are 6.42~dB, 6.54~dB, 6.67~dB for $\ell_1$ regularization, $P_{\text{EL}}$\, and $P_{\gamma}$\, respectively.  Although the SNRs are close to each other, the reconstructions show different behaviors. Both $\ell_1$ regularization and the proposed regularization have been successful in removing noise in the time-frequency regions with no activity. However since E-Lasso always keeps a component within a group, it has been less successful in suppressing noise in silent regions. 
We see that especially for higher frequencies, $\ell_1$ regularization suppresses the harmonics of the speech signal. In contrast, the proposed penalty admits a smaller weight $\lambda$ and is able to retain high frequency harmonics, while achieving a similar suppression of noise as $\ell_1$ regularization.

For $\ell_{2,1}$\, regularization and $\tilde{P}_{\gamma}$ regularization, the output SNRs are 6.09 and 7.55 dB, respectively. We observe that despite its lower SNR, $\ell_{2,1}$ preserves the harmonics better than the $\ell_1$ norm. However, especially in the low frequency region, noise suppression is modest. While this can be overcome with a higher threshold $\lambda$ for lower frequencies, such an approach requires choosing $\lambda$ in a frequency dependent manner, which complicates the application. We observe however that $\tilde{P}_{\gamma}$ achieves suppression of noise in a wide frequency range, without having to tune parameters for each frequency separately. Specifically, noise within harmonics is eliminated similarly as in the reconstructions obtained by $P_{\text{EL}}$ and $P_{\gamma}$. We also remark that specialized single-channel enhancement methods or more sophisticated penalties (such as weighted group penalties, like those in \cite{kow13p498,sie11DAFX}) can be used to achieve a superior performance than that of the proposed denoising method. Nevertheless, we think that the proposed penalty function can be used to complement or modify such alternatives as demonstrated here.

\section{Application II : Deconvolution}\label{sec:deconvolution}
In a second application, we consider a sparse deconvolution problem. In order to be able to handle an arbitrary convolution operator, we forgo convexity and consider a non-convex formulation. We provide an algorithm for the provided formulation and discuss its convergence. We also compare the performance of the penalty/threshold function with a state-of-the-art iterative thresholding method.

\subsection{A Non-Convex Formulation and a Convergent Algorithm}
Consider a minimization formulation as
\begin{equation}\label{eqn:deconvolution}
\min_x \left\{D(x) = \frac{1}{2} \| y - H\,x\|_2^2 + \lambda\,\mathbf{P}_{\gamma}(x)\right\},
\end{equation}
where $H$\, denotes a convolution operator and $\mathbf{P}_{\gamma}$ is the proposed penalty. We remark that if $H$ is not invertible, then $D(x)$ may not be convex unless $\gamma = 0$ (see the discussion in the begining of Sec.~\ref{sec:AnalysisSynthesis}). But if $\gamma = 0$, $\mathbf{P}_{\gamma}$ is the $\ell_1$ norm. For this reason, unlike  Section~\ref{sec:denoise}, we will not restrict the cost function to be convex in this section. 
We employ the forward-backward splitting algorithm (FBS) \cite{com05p168,combettes_chp, att13p91} for obtaining a local minimizer of \eqref{eqn:deconvolution}. In the current context, FBS constructs a sequence defined as,
\begin{equation}\label{eqn:ISTA}
x^{k+1} = \mathbf{T}_{(\alpha \lambda),\gamma}\bigl( x^k - \alpha H^T (Hx^k - y)\bigr).
\end{equation}
We remark that $\mathbf{T}_{(\alpha \lambda),\gamma}$ is well-defined when $\alpha\,\lambda\,\gamma < 1$. This sets an upper bound on $\alpha$. If, in addition, $\alpha < 1/\sigma(H)$, where $\sigma(H)$ denotes the spectral norm of $H$, it can also be shown using majorization-minimization techniques \cite{Lange} that the sequence in \eqref{eqn:ISTA} monotonically decreases the cost, i.e., $D(x^{k+1}) < D(x^k)$. 
Attouch et al. show in \cite{att13p91} that the algorithm \eqref{eqn:ISTA} converges under the following additional conditions,
\begin{enumerate}[(i)]
\item $P$\, is a Kurdyka-Lojasiewicz  function (\cite{att13p91}, Defn.~2.4),
\item $x^k$'s form a bounded sequence.
\end{enumerate}
 Both of these conditions are satisfied for our setup. We first remark that  the proposed penalty function $P_{\gamma}$  is continuous and for each orthant in $\mathbb{R}^n$, it can be expressed as a polynomial function. Therefore $P_{\gamma}$ (therefore $\mathbf{P}_{\gamma}$) is semi-algebraic (see Defn. 2.1 in \cite{att13p91}) and hence is a Kurdyka-Lojasiewicz function (see the discussion at the end of Sec.~2.2 in \cite{att13p91}). Also, the proposed penalty is coercive, i.e., $\mathbf{P}_{\gamma}(x)$ increases without bound as $\|x\|_2$ increases. Therefore $D(x)$ is also coercive and any sequence that monotonically decreases $D(x)$ lies in a  bounded set.  To summarize, the following proposition is a corollary of Thm.~5.1 of \cite{att13p91}.
\begin{prop}\label{prop:KL}
 If ${\alpha < \min(\sigma(H), 1/(\lambda\,\gamma))}$, then $x^k$'s in \eqref{eqn:ISTA} decrease the cost $D(x)$ monotonically, and converge to a local minimizer.\end{prop}

\subsection{Numerical Experiment}

\begin{figure}
\centering
 \includegraphics[scale=1]{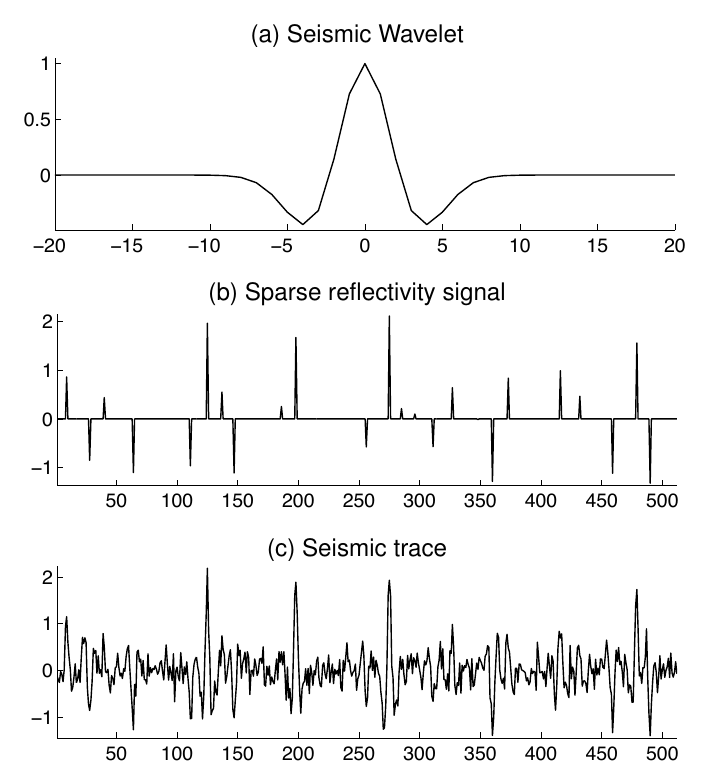}
\caption{Signals from the deconvolution experiment. (a) The seismic wavelet which is assumed to be known, (b) the sparse reflectivity signal, (c) observed noisy seismic trace. }\label{fig:SeismicObservation}
\end{figure}

In exploration seismology \cite{tak12p27,rep15p539}, the goal is to estimate an unknown reflectivity signal $x$ from the observed seismic trace $y$, which is related to $x$ as,
\begin{equation}\label{eqn:observation}
y = h\ast x + w, 
\end{equation}
where $h$\, represents the seismic wavelet and $w$ denotes white noise. Denoting the convolution operator with $h$ as $H$, we can use the formulation in \eqref{eqn:deconvolution} to estimate $x$\, from $y$. We will assume that the seismic wavelet, $h$,\, is known. Specifically, we experiment with the band-pass Ricker wavelet (dominant frequencies in the range $10 \sim 40$ Hz), sampled at $f_s = 300$~Hz, which is shown in Fig.~\ref{fig:SeismicObservation}a.  We remark that the operator $H$ used in the experiments is severely ill conditioned.

We use a synthetic sparse reflectivity signal for this experiment. The signal is selected by sampling a stochastic process where the probability of  observing a non-zero at a specific sample is 0.1, provided that a non-zero has not occured in the last 10 samples. The value of the non-zero sample is obtained by sampling a normal distribution. Notice that, this process is not a sparse Bernoulli process but a Markov process due to the  dependence on the past. The resulting synthetic $x$, of length $N = 512$ is shown in Fig.~\ref{fig:SeismicObservation}b. The observed seismic trace $y$, generated according to \eqref{eqn:observation} is shown in Fig.~\ref{fig:SeismicObservation}c. 
We used zero-mean white Gaussian noise to produce $y$. The input SNR for this observation is 5~dB.

We compared the performance of the proposed algorithm with the sparse-group lasso formulation (SGL) \cite{sim13p231} and the iterated $p$-shrinkage (IPS) algorithm \cite{woo15arXiv}, which was observed to perform very well for   sparse deconvolution (see e.g. the comparisons in \cite{sel16nonseparable}). 

SGL aims to achieve sparsity within groups and uses as few groups as possible for reconstruction. The SGL penalty is given as,
\begin{equation}
P_{\text{SGL}}(x) = \beta\,\|x\|_1 + (1-\beta)\,\sum_{i}\|x^{(i)}\|_2,
\end{equation}
where $\beta \in (0,1)$ and $\|x^{(i)}\|_2$ denotes the $\ell_2$ norm of the $i\thh$ group. Replacing $\mathbf{P}_{\gamma}$ with $P_{\text{SGL}}$ in \eqref{eqn:deconvolution}, we obtain the SGL formulation. We set $\beta = 0.95$ as in \cite{sim13p231} and make a sweep search for selecting $\lambda$. The groups consist of neighboring intervals of length $8$.  

IPS employs a threshold function depending on two parameters, namely $\lambda$ and $p$. The parameter $p$ determines the shape of the threshold function and defines a family of functions that lie between  soft ($p=1$) and hard threshold ($p\to -\infty$). We selected $p=-1/2$, which gave fairly good results. The parameter $\lambda$ is the threshold value and is selected with a sweep search.

Finally, for $\mathbf{P}_{\gamma}$, we use the same groups as SGL. We set $\gamma = 0.9/\lambda$ and select $\lambda$ with a sweep search.  We remark that for the current setup, since the distance between two non-zeros of $x$ is at least 10, each group of size 8 contains at most a single non-zero. This is the reason for choosing $\lambda\gamma$ close to unity. If multiple zeros were expected within a group, a lower value of $\gamma$ would be more feasible.

\begin{table}
\renewcommand{\arraystretch}{1.3}
\centering
\caption{SRER Performance Comparison for Deconvolution \label{table:SRER}}

\begin{tabular}{c c c c}
$\text{SNR}_{\text{in}}$& Method & $\mathbb{E}(\text{SRER})$ &  $\sigma(\text{SRER})$ \\
\hline
\multirow{3}{*}{5~dB}
& SGL & 10.08 & 0.92 \\
& IPS & 9.18 & 2.29 \\
& Proposed & 11.36 & 1.42 \\
\hline
\multirow{3}{*}{10~dB}
& SGL & 14.58 & 0.91 \\
& IPS & 14.99 & 3.36 \\
& Proposed & 16.63 & 1.51 \\
\hline
\multirow{3}{*}{15~dB}
& SGL & 19.41 & 0.87 \\
& IPS & 21.87 & 1.77 \\
& Proposed & 21.82 & 1.46 \\
\hline
\multirow{3}{*}{20~dB}
& SGL & 24.19 & 0.87 \\
& IPS & 24.08 & 1.50 \\
& Proposed & 27.09 & 1.41 \\
\hline
\end{tabular}
\end{table}

We considered four different input SNRs (5, 10, 15, 20~dB) and evaluated the deconvolution performance using signal to reconstruction error ratio (SRER), $\|x\|_2 / \| x - \hat{x}\|_2$, where $\xh$ denotes the estimate. For each input SNR value, we repeat the experiment for 500 different noise realizations to obtain average and standard deviation statistics of the performance. We set $\alpha$ to be near the upper bound allowed in Prop.~\ref{prop:KL}. We remark that the proposed formulation and IPS are essentially non-convex formulations but we have seen that both algorithms converge in our experiments (as claimed by Prop.~\ref{prop:KL} and in \cite{woo15arXiv}). We also experimented with debiasing but for the SNR range considered in these experiments, we found that debiasing actually results in lower SRER for all of the methods.

\begin{figure}
\centering
 \includegraphics[scale=1]{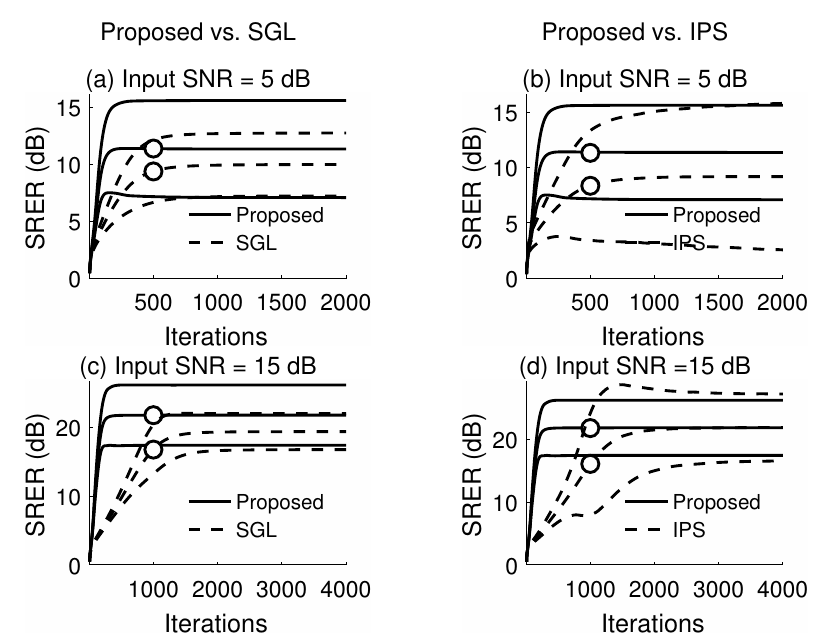}
\caption{Reliability plots comparing the proposed algorithm and iterative $p$-shrinkage (IPS) \cite{woo15arXiv}. Each figure shows the signal to reconstruction error with respect to iterations. Solid and dashed lines belong to the proposed algorithm and IPS respectively. The means are marked with white circles. The other lines indicate three times the standard deviation from the means for each method.}\label{fig:ReliabilityPlot}
\end{figure}

We observed an interesting trend with respect to different input SNRs. For low input SNRs, the proposed formulation performs better than the other two methods, in terms of average SRER. As the input SNR increases, the best SRER achieved with IPS sometimes surpasses those of the proposed formulation and SGL. However, the performance of IPS varies more with respect to different trials. On average, the proposed formulation performs better than both methods. Also, the proposed threshold function requires fewer iterations to fairly converge, although the iterations are computationally more costly. In order to visualize these, we show in Fig.~\ref{fig:ReliabilityPlot} the SRER performance with respect to iterations for input SNRs 5 and 15~dB. Here, in addition to average SRER, three times the standard deviation of the SRER with respect to iterations is also shown. For the limits, the average SRER and the standard deviation of the SRER for the different algorithms are tabulated in Table~\ref{table:SRER}. In conclusion, the proposed penalty/formulation yields a better average SRER over SGL, which targets a similar property. However, its performance is less consistent with respect to SGL as seen by a comparison of $\sigma(\text{SRER})$. On the other hand, the proposed method is favorable compared to IPS. Especially for high input SNRs, the average SRER returned by IPS can be higher but  stability is poorer. We think this is partly due to the clean sparse reflectivity signal which enforces a certain distance between the zeros. This property is taken into account by the proposed formulation but the IPS, which performs very well for arbitrary sparse signals, does not make use of this information.

\section{Conclusion}\label{sec:conc}

We proposed a group separable penalty function suitable for signals showing a sparse behavior both across and within groups of coefficients. We derived an associated threshold function $\mathbf{T}_{\lambda,\gamma}$ and studied how it behaves as its parameters vary. We argued that a good strategy is to choose $\lambda$ proportional to the noise standard deviation and $\gamma$ according to how many non-zero coefficients are expected in each group. Specifically, as $\gamma \to 1/\lambda$, we showed that in each group, there remains at most one non-zero coefficient (when the largest coefficient exceeds the threshold $\lambda$).

We think that the proposed penalty/threshold would be of interest in several areas, such as  EEG source localization \cite{gra09ISBI},  seismic deconvolution \cite{tak12p27,rep15p539}, audio processing, specifically decomposing  audio signals into transient and tonal components \cite{kow09p251,bay14p95}, 
low-rank matrix recovery \cite{can09p717,parekh15ELMA} with a bound on the rank. We hope to consider such applications in future work.

\begin{center}
\textbf{Acknowledgement}
\end{center}
We thank Pavel Rajmic, Brno University of Technology, Czech Republic, for discussions and comments.

\appendices

\section{Proof of Prop.~\ref{prop:monotone}}\label{sec:proofmonotone}
Recall that 
$\hat{x} = T_{\lambda,\gamma}(z)$ is the minimizer of $C_{\lambda,\gamma}(x|z)$ in \eqref{eqn:C} with respect to $x$.
\begin{enumerate}[(a)]
\item 
Let $z_i \geq 0$. 

Assume $\hat{x}_i > z_i$. Define a new vector $x^*$ as
\begin{align}
x^*_i &= \max\bigl(z_i -  (\xh_i - z_i) , 0 \bigr),\\
x^*_m &= \xh_m, \text{ if }m\neq i.
\end{align}
Then, $\| x^* - z\|_2^2 \leq \| \xh - z \|_2^2$ and $P(x^*) < P(\xh)$. Therefore, $C_{\lambda,\gamma}(x^*|z)< C_{\lambda,\gamma}(\xh|z)$. In words, $x^*$ achieves a strictly lower cost than $\xh$, which is a contradiction. Thus we must have, $\xh_i \leq z_i$.

Assume $\xh_i < 0$.
Define a new vector $x^*$ as
\begin{align}
x^*_i &= 0,\\
x^*_m &= \xh_m, \text{ if }m\neq i.
\end{align}
Then, $\| x^* - z\|_2^2 < \| \xh - z \|_2^2$ and $P(x^*) < P(\xh)$. Therefore, $x^*$ achieves a strictly lower cost than $\xh$, which is a contradiction. Thus we must have, $\xh_i \geq 0$. 

The second part of the claim follows similarly.
\item By part (a), we can assume without loss of generality that $z$ has non-negative components. Assume $z_i > z_m \geq 0$. Note that by part (a), $\xh_i \geq 0$, $\xh_m \geq 0$. Suppose now that $\xh_i < \xh_m$.  Define a new vector $x^*$ as
\begin{align}
x^*_i &= \xh_m,\\
x^*_m &= \xh_i,\\
x^*_l &= \xh_l, \text{ if }l\neq i \text{ or }l\neq m.
\end{align}
Then, $P_{\gamma}(x^*) = P_{\gamma}(\hat{x})$. We obtain after some algebraic manipulations that
\begin{multline}
\frac{1}{2} \left(  \|\xh - z \|_2^2 -  \| x^* - z \|_2^2 \right)  \\
 = -\left( z_m - z_i \right)\,\left( \xh_m - \xh_i \right) > 0.
\end{multline}
It thus follows $C_{\lambda,\gamma}(x^*|z) <  C_{\lambda,\gamma}(\xh|z)$,  which is a contradiction. Thus $\xh_i \geq \xh_m$.

\item Without loss of generality, assume $z$ has non-negative components. Assume also that $z_i = z_m \geq 0$. By part (a), we will have $\xh_i \geq 0$, $\xh_m \geq 0$. Suppose now that $\xh_i > \xh_m \geq 0$. Set $a = (\xh_i + \xh_m)/2$ and define a new vector $x^*$ as
\begin{align}
x^*_i &=x^*_m = a,\\
x^*_l &= \xh_l, \text{ if }l\neq i \text{ or }l\neq m.
\end{align} 
Observe that $\|\xh\|_1 = \| x^* \|_1$ and $\|\xh\|_2 > \| x^* \|_2$.
Using $z_i = z_m$, we first note,
\begin{equation}
\frac{1}{2} \left(  \|\xh - z \|_2^2 -  \| x^* - z \|_2^2 \right)  
 = \frac{1}{2} \left( \| \xh \|_2^2 - \| x^* \|_2^2 \right). 
\end{equation}
Also, since  $P_{\gamma}(\cdot) = \frac{\gamma}{2} (\| \cdot \|_1^2 - \| \cdot \|_2^2) + \| \cdot \|_1$, we find
\begin{equation}
P_{\gamma}(\xh) - P_{\gamma}(x^*) = \frac{\gamma}{2} (\|x^*\|_2^2 - \| \xh \|_2^2).
\end{equation}
Using these, we obtain, 
\begin{align}
C_{\lambda,\gamma}(\xh|z) - C_{\lambda,\gamma}(x^*|z) &= \frac{1-\lambda\gamma}{2} (\| \xh \|_2^2 - \|x^* \|_2^2) \nonumber \\ 
&> 0.
\end{align}
Thus $x^*$ achieves a lower cost than $\xh$, which is a contradiction. Therefore $\xh_i \leq \xh_m$. Changing the roles of the indices $m$ and $i$, we must also have $\xh_m \leq \xh_i$. Therefore, $\xh_i = \xh_m$.
\end{enumerate}

\section{Proof of Lemma~\ref{lem:selectk}}\label{sec:proofselectk}
\begin{enumerate}[(a)]
\item Since $z_i$'s are ordered, the assumption ${z_{i+1} > h(i)}$ implies
\begin{equation}
z_i \geq z_{i+1} > \frac{\lambda\,(1-\lambda\,\gamma) +\lambda\,\gamma\,\sum_{j=1}^i z_j  }{1 + (i-1)\,\lambda\,\gamma}.
\end{equation}
This in turn implies
\begin{equation}
z_i (1 + (i-1)\,\lambda\,\gamma ) > \lambda\,(1-\lambda\,\gamma) + \lambda\,\gamma\,\sum_{j=1}^i z_j.
\end{equation}
Subtracting $\lambda\,\gamma\,z_i$ from both sides and rearranging, we obtain
\begin{equation}
z_i >  \frac{\lambda\,(1-\lambda\,\gamma) + \lambda\,\gamma\,\sum_{j=1}^{i-1} z_j  }{1 +(i-2)\,\lambda\,\gamma} 
= h(i-1).
\end{equation}
\item The proof of this part is similar. Since $z_{i+1} \leq z_i$, the assumption $z_{i} \leq h(i)$ implies
\begin{equation}
z_{i+1} (1 + (i-1)\,\lambda\,\gamma ) \leq \lambda\,(1-\lambda\,\gamma) + \lambda\,\gamma\,\sum_{j=1}^i z_j.
\end{equation}
Adding $\lambda\,\gamma\,z_{i+1}$ to both sides and rearranging, we  obtain
\begin{equation}
z_{i+1} \leq \frac{\lambda\,(1-\lambda\,\gamma) + \lambda\,\gamma\,\sum_{j=1}^{i+1} z_j  }{1 +i\,\lambda\,\gamma} 
= h(i+1).
\end{equation}
\end{enumerate}

\section{Proof of Prop.~\ref{prop:hi}}\label{sec:proofhi}
Assume $z_{i+1} > h(i)$. This inequality implies, by the definition of $h(i)$ in \eqref{eqn:hk} that 
\begin{equation}\label{ineq1prop}
\lambda\,(1-\lambda\,\gamma) + \lambda\,\gamma\,\sum_{j=1}^i z_j  < \left(1 + (i-1)\,\lambda\,\gamma\right)\, z_{i+1}.
\end{equation}

We first show that $z_{i+1} > h(i+1)$.
Using \eqref{ineq1prop} in $h(i+1)$, we find
\begin{align}
h(i+1) &= \frac{\lambda\,(1-\lambda\,\gamma) + \lambda\,\gamma\,\sum_{j=1}^i z_j   + \lambda\,\gamma z_{i+1}}{1 + i\,\lambda\,\gamma}\\
&< \frac{\left(1 + (i-1)\,\lambda\,\gamma\right)\, z_{i+1}   + \lambda\,\gamma z_{i+1}}{1 + i\,\lambda\,\gamma} \\
& = \frac{\left(1 + i\,\lambda\,\gamma\right)\, z_{i+1} }{1 + i\,\lambda\,\gamma} \\
& = z_{i+1}.
\end{align}

Let us now show that $h(i+1) > h(i)$. Notice that for positive $a$, $b$, $c$, $d$,
\begin{equation}
\frac{a + c}{b + d} > \frac{a}{b},
\end{equation}
if and only if $ad < bc$. Now if we set
\begin{align}
a &= \lambda\,(1-\lambda\,\gamma) + \lambda\,\gamma\,\sum_{j=1}^i z_j,\\
b &= \left(1 + (i-1)\,\lambda\,\gamma\right),\\
c &= \lambda\,\gamma\,z_{i+1},\\
d &= \lambda\,\gamma,
\end{align}
then $h(i) = a/b$ and $h(i+1) = (a+c)/(b+d)$.
But we have, by \eqref{ineq1prop}
\begin{equation}
\frac{a\,d}{b\,c} = \frac{\lambda\,(1-\lambda\,\gamma) + \lambda\,\gamma\,\sum_{j=1}^i z_j }{ \left(1 + (i-1)\,\lambda\,\gamma\right)\, z_{i+1}} <1.
\end{equation}
Thus $h(i+1) > h(i)$.

\section{Proof of Prop.~\ref{prop:gamma}}\label{sec:proofgamma}

We have that ${\xh_1 \geq  \cdots \geq \xh_k >0}$ and ${\xh_{k+1}  = \cdots = \xh_n = 0}$ if and only if 
$z_{k} > h(k) > z_{k+1}$.
Plugging in the definition of $h(k)$ in \eqref{eqn:hk}, into the inequality $z_k > h(k)$, we obtain
\begin{equation}
z_k\,\bigl(1 + (k-1)\lambda\,\gamma \bigr) > \lambda\,(1-\lambda\gamma) + \lambda\gamma \sum_{i=1}^k z_i.
\end{equation}
Redistributing $z_k$'s we can write,
\begin{equation}
z_k\,(1 -\lambda\,\gamma ) > \lambda\,(1-\lambda\gamma) + \lambda\gamma \sum_{i=1}^{k-1} (z_i - z_k).
\end{equation}
This is equivalent to 
\begin{equation}
(z_k - \lambda)\,(1 -\lambda\,\gamma ) >  \lambda\gamma \sum_{i=1}^{k-1} (z_i - z_k).
\end{equation}
Dividing both sides by the positive $(z_k - \lambda) \, \lambda\,\gamma$, we find,
\begin{equation}
\frac{1}{\lambda\,\gamma} - 1 > \frac{\sum_{i=1}^{k-1} (z_i - z_k)}{z_k - \lambda}.
\end{equation}
Rearranging this equation, we obtain \eqref{eqn:5b}. 
\eqref{eqn:5a} can be shown similarly.

\section{Proof of Prop.~\ref{prop:hybrid}}\label{sec:proofhybrid}

In addition to the notation introduced in the proposition statement, let us also define $\tilde{w}$\, to be length-$m$ vector such that $\tilde{w}_i = \| \xh^{(i)}\|_2$. We first observe that if $z^{(i)} = 0$, then $\xh^{(i)} = 0$, for otherwise we could reduce the cost by setting $\xh^{(i)} = 0$. Let us denote
\begin{equation}
\partial P_{\gamma}(\tilde{w}) = ( \gamma \| \tilde{w} \|_1 + 1 )\,\sign(\tilde{w}) - \gamma \tilde{w},
\end{equation}
where $`\sign$'\, is the set valued mapping defined in \eqref{eqn:sign}.
Then, the optimality conditions for \eqref{eqn:Ttilde} imply that
\begin{equation}\label{eqn:opthybrid}
0 \in \xh^{(i)} - z^{(i)} + \lambda\,\bigl( \partial P_{\gamma}(\tilde{w}) \bigr)_i\,u^{(i)}, \text{ for }i = 1,2,\ldots, m,
\end{equation}
where $u^{(i)}$ is a unit norm vector such that $\langle \xh^{(i)}, u^{(i)} \rangle  = \|\xh^{(i)} \|_2$. That is, if $\xh^{(i)} \neq 0$, then $\xh^{(i)}  = \| \xh^{(i)}\|_2\, u^{(i)}$. This in turn implies that if $\hat{x}^{(i)} \neq 0$, then since $z^{(i)} \neq 0$ (by the observation noted above), we must also have $z^{(i)} = \|z^{(i)}\|_2 u^{(i)}$, or $\langle z^{(i)}, u^{(i)} \rangle = \| z^{(i)} \|_2$. Taking inner products with $u^{(i)}$ in \eqref{eqn:opthybrid}, we thus find,
\begin{equation}
0 \in \tilde{w}_i - w_i + \lambda\,\bigl( \partial P_{\gamma}(\tilde{w}) \bigr)_i, \text{ for }i = 1,2,\ldots, m.
\end{equation}
But these are the optimality conditions for the problem of minimizing $C_{\lambda,\gamma}(\cdot | w)$.
Therefore, $\tilde{w} = T_{\lambda,\gamma}(w)$ and the claim follows.

\end{document}